\documentclass[letterpaper, 10 pt, conference]{ieeeconf}  % Comment this line out if you need a4paper
\usepackage{multirow}
\usepackage{cite}
\usepackage{amsmath,amssymb,amsfonts} % assumes amsmath package installed
\usepackage[linktocpage=true,colorlinks=true,linkcolor=blue,citecolor=blue,urlcolor=blue]{hyperref}
\usepackage{tikz}
\usepackage{ifthen}

\usepackage{fancyvrb}
\usetikzlibrary{shapes}
\usetikzlibrary{calc}

\providecommand{\keywords}[1]
{
  \small	
  \textbf{\textit{Keywords---}} #1
}

\newcommand{\bx}{\boldsymbol{x}}

\newcommand{\blambda}{\boldsymbol{\lambda}}

\newcommand{\bu}{\boldsymbol{u}}
\newcommand{\bw}{\boldsymbol{w}}

\newcommand{\bz}{\boldsymbol{z}}

\newcommand{\bZ}{\boldsymbol{Z}}
\newcommand{\bv}{\boldsymbol{v}}

\newcommand{\bS}{\boldsymbol{S}}

\newcommand{\bd}{\boldsymbol{d}}
\newcommand{\bp}{\boldsymbol{p}}
\newcommand{\bs}{\boldsymbol{s}}

\newcommand{\bh}{\boldsymbol{h}}

\newcommand{\bA}{\boldsymbol{A}}

\newcommand{\bB}{\boldsymbol{B}}

\newcommand{\bH}{\boldsymbol{H}}

\newcommand{\bJ}{\boldsymbol{J}}

\newcommand{\br}{\boldsymbol{r}}

\newcommand{\bSigma}{\boldsymbol{\Sigma}}

\newcommand{\diag}{\mathop{\text{\normalfont diag}}}

\usepackage{algorithm}
\usepackage{algpseudocode}

\algdef{SE}[SUBALG]{Indent}{EndIndent}{}{\algorithmicend\ }%
\algtext*{Indent}
\algtext*{EndIndent}

\usepackage[utf8]{inputenc} % allow utf-8 input
\usepackage[T1]{fontenc}    % use 8-bit T1 fonts
\usepackage{hyperref}       % hyperlinks
\usepackage{url}            % simple URL typesetting
\usepackage{booktabs}       % professional-quality tables
\usepackage{amsfonts}       % blackboard math symbols
\usepackage{nicefrac}       % compact symbols for 1/2, etc.
\usepackage{microtype}      % microtypography

\IEEEoverridecommandlockouts                              % This command is only needed if
\title{\LARGE \bf
  Exploiting GPU/SIMD Architectures for Solving\\ Linear-Quadratic MPC Problems*
  \thanks{*This material is based upon work supported by the U.S. Department of Energy, Office of Science, Office of Advanced Scientific Computing Research (ASCR) under Contract DE-AC02-06CH11347.}% <-this % stops a space
}

\author{David Cole,$^{1}$ Sungho Shin,$^{2}$ François Pacaud,$^{2}$ Victor M. Zavala,$^{1,2}$ Mihai Anitescu$^{2,3}$% <-this % stops a space
  \thanks{$^{1}$Department of Chemical and Biological Engineering, University of Wisconsin-Madison, Madison, WI 53706}%
  \thanks{$^{2}$Mathematics and Computer Science Division, Argonne National Laboratory, Lemont, IL 60439}%
  \thanks{$^{3}$Department of Statistics, University of Chicago, Chicago, IL 60637}
}

\begin{document}

\maketitle
\thispagestyle{empty}
\pagestyle{empty}

%%%%%%%%%%%%%%%%%%%%%%%%%%%%%%%%%%%%%%%%%%%%%%%%%%%%%%%%%%%%%%%%%%%%%%%%%%%%%%%%
\begin{abstract}
  We report numerical results on solving constrained linear-quadratic model predictive control (MPC) problems by exploiting graphics processing units (GPUs).
  The presented method reduces the MPC problem by eliminating the state variables
  and applies a condensed-space interior-point method to remove the inequality
  constraints in the KKT system. The final condensed matrix is positive
  definite and can be efficiently factorized in parallel on GPU/SIMD architectures.
  In addition, the size of the condensed matrix depends only on the number of controls in the
  problem, rendering the method particularly effective when the problem has many states but few inputs and moderate horizon length.
  Our numerical results for PDE-constrained problems show that the approach is an order of magnitude faster than a standard CPU implementation.
  We also provide an open-source Julia framework that facilitates  modeling ({\tt DynamicNLPModels.jl}) and solution ({\tt MadNLP.jl}) of MPC problems on GPUs. 
\end{abstract}

%%%%%%%%%%%%%%%%%%%%%%%%%%%%%%%%%%%%%%%%%%%%%%%%%%%%%%%%%%%%%%%%%%%%%%%%%%%%%%%%
\section{INTRODUCTION}

The fundamental challenge of model predictive control (MPC) is solving optimal control (dynamic optimization) problems within short sample time intervals. The real-time computation load can be prohibitive, especially when the system dimension is high or the prediction horizon is long; such computational challenges have limited the application scope of MPC. Accordingly, the scalable solution of optimal control problems has been a long-standing challenge in MPC, with gradual improvements being made to the algorithm \cite{wright1991partitioned,rao1998application,frison2014high,lee2022gpu,andersson2013condensing,jerez2012sparse,na2022convergence}.

While numerical solvers have greatly improved thanks to strides in algorithms and computing hardware, the performance of single-core processors has started to stall in the past decade.
Rather, progress in hardware has been primarily driven by parallel architectures such as multicore processors, distributed computing clusters, and graphics processing units (GPUs). Hence, in order to continue leveraging advances in modern computing hardware,  algorithms and software implementations that harness such capabilities must be developed. While the use of multicore processors and distributed computing has been widely studied in the field of mathematical optimization \cite{schnabel1995view,wright2001solving,chiang2014structured,kim2019asynchronous},
relatively fewer contributions have been made regarding the use of GPUs and, more specifically,  the use of single instruction, multiple-data (SIMD) architectures. 
%As these parallel architectures are at the forefront of computational advancements, we seek to make optimal control problems more readily solvable on GPU/SIMD architectures. Doing so could make these architectures more practical in other areas, such as embedded systems (e.g., robotics and autonomous vehicles), where the system requires an efficient, on-board controller that makes decisions. Furthermore, GPU/SIMD architectures have the potential to decrease energy use, as compared to CPUs \cite{Dong2014hydrodynamic,Qasaimeh2019comparing}.
In this context, our work is devoted to accelerating the solution of linear-quadratic MPC problems by exploiting modern GPU/SIMD architectures. 

Developing algorithms for MPC problems on GPU/SIMD architectures offers another benefit for embedded applications (e.g., robotics and autonomous vehicles), in addition to reducing the solution time.
The embedded systems with GPU/SIMD capabilities, such as NVIDIA Jetson or the BeagleBoard X15, have the potential to decrease energy use, as compared with CPUs \cite{Qasaimeh2019comparing}.
Accordingly, they 
allow the potential to  create new applications for embedded systems \cite{Dong2014hydrodynamic,forgione2020efficient}.
% Thus, these architectures could greatly improve embedded systems (e.g., robotics and autonomous vehicles) where the system requires an efficient, on-board controller that makes decisions \cite{}.
% developing tools for the efficient solution of linear-quadratic MPC on GPU could.
% These can be implemented on the systems like the .
The use of such embedded systems for MPC has been reported in recent works \cite{phung2017model, rathai2020gpu, yu2017efficient}.
% Also, several studies have explored using the parallel architecture of field programmable gate arrays and shown promising results \cite{jerez2011parallel, ling2008embedded, lucia2017optimized}. 

% Thus, accelerating the solution of linear-quadratic MPC through GPU/SIMD architectures could open new applications for embedded systems.

The efficient solution of MPC problems (on either CPUs or GPUs) requires tailored linear algebra techniques that exploit the structure of the Karush--Kuhn--Tucker (KKT) systems arising in the optimization procedure.
Classical linear algebra methods include Riccati-like recursions \cite{frison2014high}, sparse LU factorizations \cite{lee2022gpu}, or sparse LDL$^\top$ factorizations \cite{HSL}.
These methods rely on direct linear solvers that exploit the sparsity or the block tridiagonal structure of the KKT systems. Despite being highly efficient on CPUs, most sparse matrix factorization routines are sequential and known to be difficult to parallelize. Furthermore, when the factorization of the sparse matrix requires a lot of fill-ins (e.g., as in PDE systems \cite{schenk2009inertia}), the sparse direct solver can become extremely slow.

An alternative approach to sparse methods is to formulate the problem in a reduced space by eliminating the state variables. This is readily done by writing each state as a function of the initial state and the previous controls \cite{jerez2012sparse} (\emph{reduction}). By further eliminating the inequality constraints using a Schur complement technique (\emph{condensation}), the Newton step computation can be performed by solving a small (the only variables are controls) positive definite system, which can be factorized efficiently with dense Cholesky. Jerez et al.~\cite{jerez2012sparse} highlight this method and compare the computational complexity of the step computation within sparse and dense methods. Their complexity analysis suggests that the dense method is particularly effective when the number of states is large, the number of controls is small, and the time horizon is short. We note that the dense formulation has been used in various contexts \cite{rawlings2017, frison2016algorithms, andersson2013condensing} and is implemented in the state-of-the-art MPC solver HPIPM \cite{frison2020hpipm}.

All these methods transpose directly in the linear algebra routine employed inside the optimization solver
and can be exploited in conjunction with a GPU-accelerated interior-point method (IPM) \cite{jung2008implementing, smith2011gpu, cao2016augmented, shah2020efficiency, pacaud2022condensed, lee2022gpu}. Indeed, in direct contrast to active-set methods (which imply expensive reordering operations in the KKT systems, associated with changes in the active set), IPM involves the solution of a sequence of KKT systems with a fixed sparsity pattern. Hence, the computational burden lies primarily at the linear algebra level (i.e., parallelizing IPM requires solving KKT systems in parallel). To overcome the limitations of sparse direct indefinite factorization routines  on GPUs \cite{swirydowicz2022linear},
 attempts to solve IPM on GPUs have relied on alternative linear algebra routines such as iterative methods \cite{smith2011gpu}, CPU-GPU methods \cite{jung2008implementing}, or reduction methods \cite{pacaud2022condensed}.
In the context of MPC, Gade-Nielsen and co-workers \cite{gade2012mpc, gade2014interior} developed an algorithm and implementation for solving linear programs for MPC (LP-MPC) by factorizing the matrix on a GPU.
In their case study, they reported an order of magnitude speedup on the GPU as their problem scale increased. Lee et al.~\cite{lee2022gpu} applied an IPM to the standard linear-quadratic MPC formulation with sparse LU factorizations on the GPU. Their method resulted in speedups of a factor of 3--4 over an alternative implementation on the CPU.

A different approach for implementing optimization algorithms on  GPU/SIMD architectures is to use iterative linear solvers. Recently, a new algorithm for general quadratic programs (QPs) based on the operator-splitting technique has been proposed  \cite{stellato2020osqp}.
The technique eliminates the necessity to factorize the KKT system in each iteration and greatly accelerates the solution.
% however, achieving high solution accuracy with iterative linear solvers is challenging.
% This method is known to be one of the fastest general-purpose algorithms that are currently available for convex QP.
The algorithm was subsequently adapted for GPUs by replacing the initial sparse factorization step (difficult to parallelize) with an iterative method (in particular, preconditioned conjugate gradient method) \cite{osqp-gpu}.
It has been reported that more than $10$ times speedup can be made through GPU accelerations \cite{osqp-gpu}.

Our goal is to demonstrate the capability of modern GPU architectures for solving dense linear-quadratic MPC problems. We formulate the problem by applying the reduction method presented in \cite{jerez2012sparse} and applying the condensed-space interior-point method with dense Cholesky factorization routines within the CUDA library.\footnote{
  % To avoid potential confusion in the terminology,
  We highlight that state elimination (\emph{reduction to a dense formulation}) is performed before solving the problem and that the inequality elimination (\emph{condensation}) is performed internally within the interior-point solver when formulating the linear system for the step computations.
  % When we say ``dense formulation'', we always mean by the reduced MPC problem whose state variables are eliminated, and when we say ``condensation'', we always mean by the inequality elimination procedure.
}
To the best of our knowledge, this is the first implementation of condensed-space IPM running on GPU for solving linear-quadratic MPC problems and with the potential to be extended to nonlinear MPC.
The \emph{dense problem formulation} is implemented in the modeling library {\tt DynamicNLPModels.jl}, and the \emph{condensation} is implemented in the nonlinear optimization solver {\tt MadNLP.jl}.
The proposed method is demonstrated by using case studies for PDE-constrained optimal control problems.
Our results show that solving dense problems using the GPU can offer an order of magnitude speedup.

\section{Dense Formulation}\label{sec:condensed_form}

We consider the following linear-quadratic MPC problem:
\begin{subequations}\label{eqn:sparse}
  \begin{align}
    \min_{\substack{\{x_t\}_{t=0}^T, \\ \{u_t\}_{t=0}^{T-1}}}
    &\; x_T^\top Q_f x_T + \sum_{t = 0}^{T-1}
      \begin{bmatrix} x_t \\ u_t \end{bmatrix} ^\top
    \begin{bmatrix} Q & S \\ S^\top & R \end{bmatrix}\begin{bmatrix} x_t \\ u_t \end{bmatrix}\\
    \textrm{s.t.}
    &\; x_0 = \overline{x} ,\\
    &\;x_{t+1} = Ax_t + Bu_t + w_t \quad \forall t \in \mathbb{I}_{[0,T-1]} \label{eqn:dyn}\\
    % &\;  \\
    &\; g^l \le E x_t + F u_t \le g^u \quad \forall t \in \mathbb{I}_{[0,T-1]} \label{eqn:ineq_constraints}\\
    &\; x^l \le x_t \le x^u \quad \forall t \in \mathbb{I}_{[0,T]} \label{eqn:state_bounds}\\
    &\; u^l \le u_t \le u^u \quad \forall t \in \mathbb{I}_{[0,T-1]}, \label{eqn:input_bounds}
  \end{align}
\end{subequations}
where $x_t\in\mathbb{R}^{n_x}$ are the states at time $t$, $u_t\in\mathbb{R}^{n_u}$ are the inputs at time $t$, $\overline{x}$ is the initial known state, and $T$ is the number of steps in the time horizon. Here $x^l,x^u,w_t \in \mathbb{R}^{n_x}$, $u^l,u^u \in \mathbb{R}^{n_u}$, $g^l,g^u \in \mathbb{R}^{n_c}$, $A \in \mathbb{R}^{n_x \times n_x}$, $B \in \mathbb{R}^{n_x \times n_u}$, $Q ,Q_f \in \mathbb{R}^{n_x \times n_x}$, $R \in \mathbb{R}^{n_u \times n_u}$, $S \in \mathbb{R}^{n_x \times n_u}$, $E \in \mathbb{R}^{n_c \times n_x}$, $F \in \mathbb{R}^{n_c \times n_u}$ are problem data. We use $n_x$ for the number of states, $n_u$ for the number of inputs, and $n_c$ for the number of constraints in Equation \eqref{eqn:ineq_constraints}.

The dense form of \eqref{eqn:sparse} can be derived by eliminating the state variables. The detailed procedure can be found in \cite{jerez2012sparse};  here we  outline the key ideas. First, to allow compact notation, we let $\bx := [x_0^\top, x_1^\top, \cdots, x_T^\top]^\top$, $\bu := [u_0^\top, u_1^\top, \cdots, u_{T-1}^\top]^\top$, and $\bw := [w_0^\top, w_1^\top, \cdots, w_{T-1}^\top]^\top$. For inducing the sparsity (when $(A,B)$ is controllable) and for the numerical stability in the elimination procedure (note that the power of $A$ diverges quickly if $A$ is unstable), we introduce a new variable $v_t$ with stabilizing feedback $K$ by
\begin{align}\label{eqn:K_constraint}
  u_t = Kx_t + v_t \quad  \forall t \in \mathbb{I}_{[0,T-1]},
\end{align}
and we let $A_K := A + BK$ and $\bv := [v_0^\top, v_1^\top, \cdots, v_{T-1}^\top]$. After the elimination of the state variables, the dynamic constraint \eqref{eqn:dyn} reduces to
\begin{align}\label{eqn:x_def}
    \bx = \bA \overline{x} + \bB \bv + \tilde{\bA} \bw,
\end{align}
where $\bA$, $\tilde{\bA}$, and $\bB $ are defined as
\begin{align*}
  \bA &:= \begin{bmatrix}
    I_n \\ A_K \\ A_K^2 \\ \vdots \\ A_K^{T-1} \\ A_K^T
  \end{bmatrix},
  \tilde{\bA} :=
  \begin{bmatrix}
     &   &  &  &  \\
    I_n &   &  &  &  \\
    A_K & I_n  & &  & \\
    \vdots &  & \ddots & & \\
    A_K^{T-2} & A_K^{T-3}  & &  I_n &  \\
    A_K^{T-1} & A_K^{T-2}  &\cdots& A_K & I_n
  \end{bmatrix}\\
  \bB &:=  \begin{bmatrix}
     &   & &  &  \\
    B &   &  &  &  \\
    A_K B & B  & &  & \\
    \vdots &  & \ddots & & \\
    A_K^{T-2} B & A_K^{T-3}B  & &  B &  \\
    A_K^{T-1} B & A_K^{T-2}B  &\cdots& A_K B & B
  \end{bmatrix}.
\end{align*}
By substituting \eqref{eqn:K_constraint} and \eqref{eqn:x_def} into \eqref{eqn:sparse}, we obtain the dense MPC problem:
\begin{subequations}\label{eqn:dense}
    \begin{align}
        \min_{\bv} &\;\; \frac{1}{2} \bv^\top \bH \bv + \bh^\top \bv + \bh_0\\
        \textrm{s.t.} &\; \bJ \bv \leq \bd,\label{eqn:dense-ineq}
    \end{align}
\end{subequations}
where $\bH$, $\bJ$, $\bh$, $\bh_0$, and $\bd$ are defined implicitly by the elimination procedure.

Observe that the problem is expressed only in terms of $\bv$. This dense formulation has $Tn_u$ variables and $2T(n_c+n_x+n_u)$ inequality constraints. Thus, the number of variables is significantly reduced, but now the problem has dense Hessian $\bH$ and Jacobian $\bJ$. By exploiting the controllability, one can partially maintain the sparsity; for more details on the state elimination procedure, readers are pointed to \cite{jerez2012sparse}.

\section{Condensed-Space Interior-Point Method}\label{sec:IPM}

We now describe the condensed-space IPM applied to the dense MPC problem in \eqref{eqn:dense}. First, we reformulate \eqref{eqn:dense} by introducing the slack variable $\bs$:
\begin{subequations}
    \begin{align}
        \min_{\bv, s} &\; \frac{1}{2}\bv^\top \bH \bv + \bh^\top \bv + \bh_0\\
        \textrm{s.t.} &\; \bJ \bv - \bd + \bs = 0 \label{eqn:slack_constraints}\\
        &\; \bs \ge 0.\label{eqn:slack}
    \end{align}
\end{subequations}
%where $\bs$ is the slack variable. 
The IPM treats the inequality constraints by replacing the inequality \eqref{eqn:slack} by the log barrier function $-\mu \sum_{i = 1}^m \log(s_i)$, where $\mu$ is the barrier term and $s_i$ is the $i$th value of $\bs$. The resulting barrier subproblem has the form
\begin{subequations}\label{eqn:barrier_subproblem}
    \begin{align}
      \min_{\bv, s} &\; \frac{1}{2} \bv^\top \bH \bv + \bh^\top \bv + \bh_0 - \mu \sum_{i = 1}^m \log(s_i) \\
      \textrm{s.t.} &\; \bJ \bv - \bd + \bs = 0\label{eqn:barrier_constraints}.
    \end{align}
\end{subequations}
The resulting KKT conditions are then
\begin{subequations}\label{eqn:KKT_conditions}
    \begin{align}
      \bH\bv + \bh +  \bJ^\top \blambda= 0 \label{eqn:KKT_cond1}\\
      \blambda - \bz  = 0\label{eqn:KKT_cond2}\\
      \bJ \bv - \bd + \bs = 0\label{eqn:KKT_cond3},
    \end{align}
\end{subequations} % need to check again
where $\bz:= \mu \bS^{-1} \boldsymbol{1}$, $\blambda$ is the Lagrange multiplier of \eqref{eqn:barrier_constraints}, $\bS := \diag(\bs)$, and $\boldsymbol{1}$ is the vector of ones.

Typical IPMs compute the step direction by applying Newton's method to the primal-dual equation in \eqref{eqn:KKT_conditions}. The step computation involves the solution of linear systems of the following form:
\begin{equation}\label{eqn:aug_KKT}
  \begin{bmatrix}
    \bH&  & \bJ^\top \\
     & \bSigma_s & I \\
    \bJ & I &
  \end{bmatrix}
  \begin{bmatrix}
    \bp^v \\ \bp^s \\ \bp^\lambda
  \end{bmatrix}  = -
  \begin{bmatrix}
    \br_1 \\ \br_2 \\ \br_3 ,
  \end{bmatrix}
\end{equation}
where $\bp^v$, $\bp^s$, and $\bp^\lambda$ are the descent directions for $\bv$, $\bs$, and $\bz$, respectively; $\br_1 := \bH\bv + \bh + \bJ^\top \blambda$, $\br_2 := \blambda - \mu \bS^{-1}\boldsymbol{1}$, $\br_3 := \bJ \bv - \bd + \bs$; and $\bSigma_s := \bS^{-1} \bZ$ with $\bZ := \diag(\bz)$. Note that this requires that $\bS$ be invertible, but this is always the case in IPMs because the primal iterate always stays in the strict interior $\bs > 0$.

This KKT system can be further condensed by writing $\bp^\lambda$ as a function of $\bp^s$ and $\bp^s$ as a function of $\bp^v$. We then obtain a system of equations for $\bp^v$:
\begin{equation}\label{eqn:cond_KKT}
    \left(\bH + \bJ^\top \bSigma_s \bJ \right) \bp^v = -\br_1 + \bJ \br_2 - \bJ^\top \bSigma_s \br_3.
\end{equation}
Here, $\bH + \bJ^\top \bSigma_s \bJ$ is positive definite and can be factorized with the Cholesky factorization, and $\bp^v$ can be computed by back triangular solve. Once $\bp^v$ is obtained, $\bp^s$ can be obtained from $\bp^s = -\br_3 - \bJ \bp^v$, and $\bp^\lambda$ can be found from $\bp^\lambda = -\br_2 + \bSigma_s \br_3 + \bSigma_s \bJ \bp^v$. After the descent directions have been determined, a backtracking line search is performed to determine the step size, and the next iterate is obtained based on the step size and the direction. The new values of $\bv$, $\bs$, and $\bz$ at the end of this process are then used in the next iteration, and the process is continued until a stopping criterion is satisfied. The method is summarized in Algorithm \ref{al:condensed_IPM}; here, the superscript $(k)$ represents the iterate at $k$th iteration.

\begin{algorithm}[t]
  \begin{algorithmic}
    \Require $A,B,Q,R,S,K,N,x_0$
    \State{Construct} $\bH$, $\bh$, $\bh_0$, $\bJ$, and $\bd$
    \State{Initialize} $\bv^{(0)}$, $\bs^{(0)}$, and $\blambda^{(0)}$
    \While{Termination criteria are not met}
    \State ${\bSigma_s}^{(k)} \leftarrow (\bS^{(k)})^{-1}\bZ^{(k)}$
    \State Compute $\bH + \bJ^\top {\bSigma_s}^{(k)} \bJ$
    \State Factorize $\bH + \bJ^\top {\bSigma_s}^{(k)} \bJ$ and solve \eqref{eqn:cond_KKT} for $\bp^{v(k)}$.
    \State $\bp^{s(k)} \leftarrow -\br_3 - \bJ \bp^{v(k)}$
    \State $\bp^{\lambda(k)} \leftarrow -\br_2 + \bSigma_s \br_3 + \bSigma_s \bJ \bp^{v(k)}$
    \State $\bp^{z(k)}\leftarrow \mu \bS^{-1}\boldsymbol{1} - \bz^{(k)} - \bSigma \bp^{s(k)}$
    \State Determine the step length $\alpha$, $\alpha_z$ through line search
    \State $\bv^{(k)}\leftarrow \bv^{(k)} + \alpha \bp^{v(k)}$\\
    \State $\bs^{(k)}\leftarrow \bs^{(k)} + \alpha \bp^{s(k)}$\\
    \State $\blambda^{(k)}\leftarrow \blambda^{(k)} + \alpha \bp^{\lambda(k)}$
    \State $\bz^{(k)}\leftarrow \bz^{(k)} + \alpha_z\bp^{z(k)}$
    \EndWhile
  \end{algorithmic}
  \caption{Condensed-space interior-point method}
  \label{al:condensed_IPM}
\end{algorithm}

\begin{figure}[t]
    \centering
    \includegraphics[width=.5\textwidth]{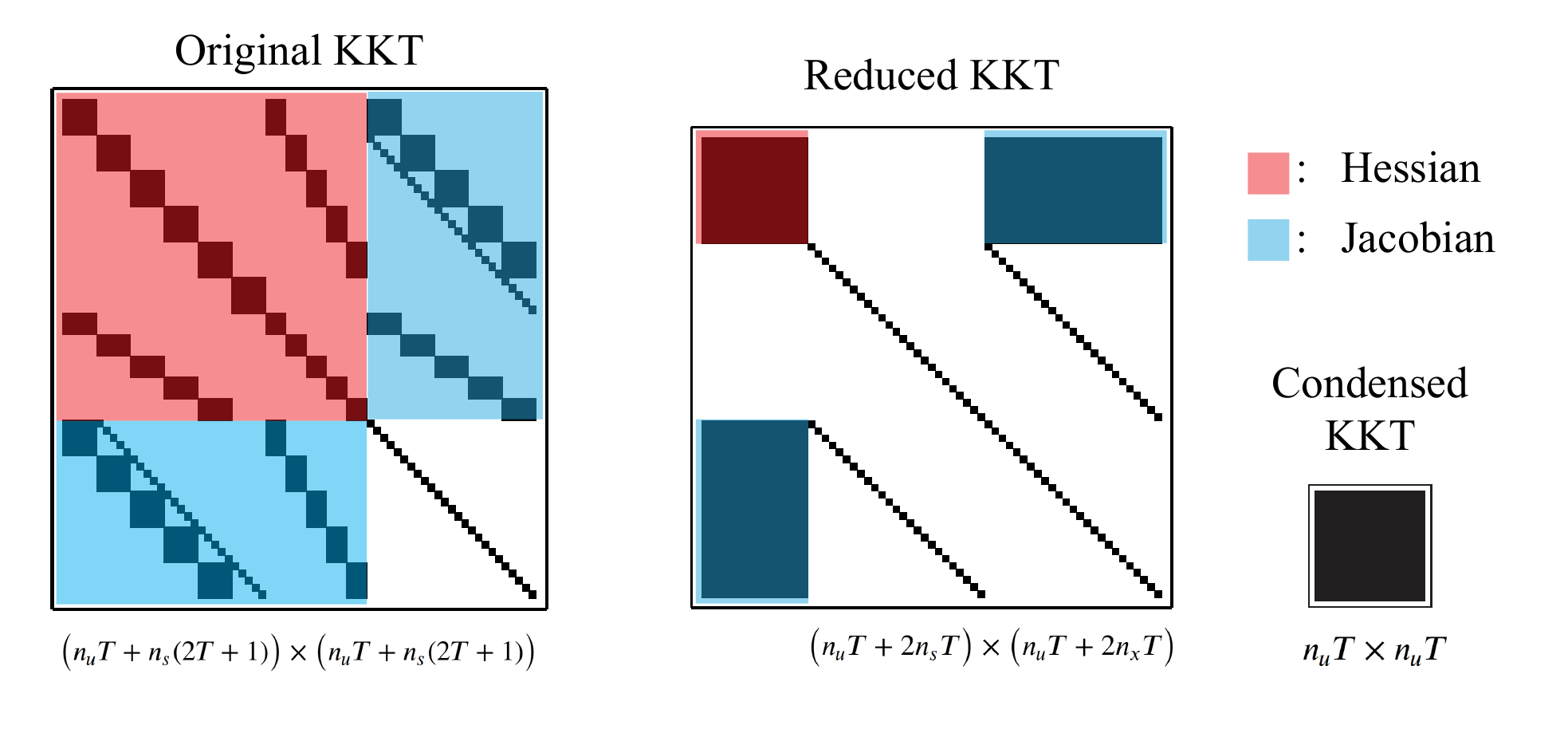}
    \caption{Full, reduced, and condensed system for linear-quadratic MPC problems.}
    \label{fig:kkt_comp}
\end{figure}

The system in \eqref{eqn:cond_KKT} is (i) smaller ($n_uT$ rows and columns) than the system in \eqref{eqn:aug_KKT} ($n_uT + 2(n_c + n_x + n_u)T$ rows and columns), (ii) dense, and (iii) positive definite. These features make it more efficient to perform Algorithm \ref{al:condensed_IPM} on the GPU or other parallel architecture because the Cholesky factorization of small dense systems can be easily parallelized, and efficient implementation is available within the CUDA library. An illustration of the sparsity pattern of full, reduced, and condensed KKT systems with one-sided bounds on the state variables is shown in Figure \ref{fig:kkt_comp} (for $n_x = 5$, $n_u=3$, $n_c = 0$, and $T = 5$).

\section{Implementation}

To facilitate building and solving sparse or reduced linear-quadratic MPC problems, we have developed the Julia package {\tt DynamicNLPModels.jl}. This package takes data given in \eqref{eqn:sparse} and builds either the sparse or reduced forms of these problems. The package is open source, and the source code can be found at \cite{dnlp}. Further, the model returned by {\tt DynamicNLPModels.jl} is a subtype of the generic NLP model data structure implemented in {\tt NLPModels.jl} \cite{nlpmodels2020} and so is compatible with other solvers in the {\tt JuliaSmoothOptimizer} ecosystem. By building the reduced problem outside the solver, {\tt DynamicNLPModels.jl} facilitates forming the reduced KKT system \eqref{eqn:aug_KKT}, which can then be further condensed \eqref{eqn:cond_KKT} as the system in \eqref{eqn:cond_KKT} within {\tt MadNLP.jl}.

{\tt DynamicNLPModels.jl} allows the user to easily define either the sparse (Equation \eqref{eqn:sparse}) or reduced (Equation \eqref{eqn:dense}) linear-quadratic MPC problem by passing user-defined data. An example of how these models are built can be seen in the code snippet below.

\begin{small}
\begin{Verbatim}[frame=single]
using DynamicNLPModels
using Random, LinearAlgebra

Q  = 1.5 * Matrix(I, (3, 3))
R  = 2.0 * Matrix(I, (2, 2))
A  = rand(3, 3)
B  = rand(3, 2)
T  = 5
s0 = [1.0, 2.0, 3.0]
                
sparse_lqdm = 
  SparseLQDynamicModel(s0, A, B, Q, R, T)
dense_lqdm  = 
  DenseLQDynamicModel(s0, A, B, Q, R, T)
\end{Verbatim}
\end{small}
The sparse and reduced models can be constructed by passing the user-defined data to {\tt SparseLQDynamicModel} or {\tt DenseLQDynamicModel}. For the dynamic model constructors, keyword arguments can also be passed for the data given in Equations \eqref{eqn:sparse} and \eqref{eqn:K_constraint}, including $Q_f$, $S$, $K$, $E$, $F$, $x^l$, $x^u$, $u^l$, $u^u$, $\bw$, $g^l$, and $g^u$. 

% \begin{figure}\small
%   \label{fig:code_snippet}
% \begin{verbatim}
% using DynamicNLPModels
% using Random, LinearAlgebra

% Q  = 1.5 * Matrix(I, (3, 3))
% R  = 2.0 * Matrix(I, (2, 2))
% A  = rand(3, 3)
% B  = rand(3, 2)
% T  = 5
% s0 = [1.0, 2.0, 3.0]
                
% lqdd = LQDynamicData(s0, A, B, Q, R, T)
                
% sparse_lqdm = SparseLQDynamicModel(lqdd)
% dense_lqdm  = DenseLQDynamicModel(lqdd)
        
% # or 
        
% sparse_lqdm = 
%   SparseLQDynamicModel(s0, A, B, Q, R, T)
% dense_lqdm  = 
%   DenseLQDynamicModel(s0, A, B, Q, R, T)
% \end{verbatim}
%   \caption{Code snippet for building a sparse or reduced (dense) linear-quadratic MPC problem using {\tt DynamicNLPModels.jl}}
% \end{figure}

\section{Numerical Results}

We apply the condensed-space IPMs to a 3-D PDE temperature control problem to showcase the potential benefits and limitations of our approach. For problems with low numbers of inputs and high numbers of constraints, we show that solving the dense linear-quadratic MPC problem on the GPU using a condensed-space IPM can offer an order of magnitude speedup over solving the dense problem on the CPU with the condensed-space IPM or solving the sparse problem on the CPU. Also, our results suggest that the condensed-space GPU approach is less effective for problems with a low ratio of states to inputs. In all cases, we constructed the model using {\tt DynamicNLPModels.jl}. We solved the sparse form of the problem using {\tt Ipopt} \cite{wachter2006implementation} with MA27 \cite{HSL} as the linear solver. We solved the dense form of the problem on the CPU as well as the GPU using {\tt MadNLP.jl} and dense Cholesky factorization using LAPACK (OpenBLAS) or CUSOLVER. Scripts and solver outputs are available at \cite{code}.

\subsection{Problem Formulation}

We adapted the 2-D thin plate temperature control problem given by \cite{na2022convergence} to a 3-D temperature control problem in a cube using the boundary inputs. In this problem, each face of the cube can be heated by control to try to reach the desired temperature profile. The continuous form of this problem can be stated as follows:

\begin{subequations}\label{eqn:3DPDE}
    \begin{align}
      \small \min_{x, u} &\small \; \frac{1}{2} \int_0^{T} \int_{w \in \Omega} \Bigg( q(x(w, t) - d(w, t))^2  dw  \\
                         &\qquad\qquad\qquad\qquad+ \sum_{i \in \mathcal{U}} r u_i(t)^2 \Bigg) dt\nonumber\\
        \begin{split}
        \textrm{s.t.} &\; \rho C_p  \frac{\partial x(w, t)}{dt} = k \nabla^2 x(w, t), \; w \in \Omega, t \in [0, T]
        \end{split}\\
        &\; x^l \le x(w, t) \le x^u \; w \in \Omega,\; t \in [0, T]\\
        &\; u^l \le u_i(t) \le u^u \; i \in \mathcal{U},\; w \in \Omega, \;t \in [0, T]\\
        &\; x({w}, t) = u_i(t), \; i \in \mathcal{U},\; {w} \in \partial \Omega_i, \; t\in[0,T] \label{eqn:dirichlet_cond} \\
        &\; x(w, 0) = \hat{x}, \; w \in \Omega,
    \end{align}
  \end{subequations}
  where $\Omega = [0, L]^3 \subseteq \mathbb{R}^3$ is the 3-D domain of a cube of length $L$, $\mathcal{U}:=\{1,2,\cdots,6\}$ is the set of all cube faces, $x: \Omega \times [0, T] \rightarrow \mathbb{R}$ is the temperature, $d: \Omega \times [0, T] \rightarrow \mathbb{R}$  is the set point temperature,  $u_i: [0, T] \rightarrow \mathbb{R}$ are the boundary temperatures for $i \in \mathcal{U}$, and $\partial \Omega_i \subseteq \Omega$ for $i\in\mathcal{U}$ are the boundaries. Note that \eqref{eqn:dirichlet_cond} enforces the Dirichlet boundary condition. The constant values are given in Table \ref{tab:const}.

  \begin{table}[t]
    \centering
    \caption{Parameter values for the 3-D temperature control problem.}
    \label{tab:const}
    \begin{tabular}{llr}
      \hline\hline
      $\Delta_t$ &Temporal discretization mesh size & 0.1 sec\\
      \hline
      $\Delta_w$ &Spatial discretization mesh size & 0.02 m\\
      \hline
      %r, q, ell_x, xl, xu, ul, uu, rho, cp, k, xhat
      $L$ & Length of cube face & $(N + 1) \Delta_w$ \\
      \hline
      $\rho$ & Density of copper & 8960 kg $\cdot$ m\textsuperscript{-3}\\
      \hline
      $C_p$ & Specific heat of copper & 386 J $\cdot$ (kg $\cdot$ K)\textsuperscript{-1}\\
      \hline
      $k$ & Thermal conductivity of copper & 400 W $\cdot$ (m $\cdot$ K)\textsuperscript{-1}\\
      \hline
      $q$& Weight on temperature error &10 $\Delta_w^2 $ m\textsuperscript{-2}\\
      \hline
      $r$ & Weight on input & $\frac{1}{10} \Delta_w^2 $ m\textsuperscript{-2}\\
      \hline
      $(x^l, x^u)$ &(Lower, upper) temperature bounds & $(200, 550)$ K\\
      \hline
      $(u^l, u^u)$ & (Lower, upper) input bounds & $(300, 500)$ K\\
      \hline
      $\hat{x}$ & Initial temperature & $300$ K\\
      \hline\hline
    \end{tabular}
  \end{table}

  The infinite-dimensional problem in \eqref{eqn:3DPDE} is discretized by using a finite difference method in space and an explicit Euler scheme in time to make a finite-dimensional MPC problem of the form in \eqref{eqn:sparse}.
  We also treat the boundaries as inputs only (not as states), so all of the state variables are strictly internal temperatures.
  In addition, this problem can be scaled up by changing $T$ (the time horizon) or $N$ (the number of discretization points per dimension).
  Note that increasing $N$ increases cubically the number of states because the spatial domain is three-dimensional.

\subsection{Results}
\begin{figure}[t]
    \centering
    \includegraphics[scale=.22]{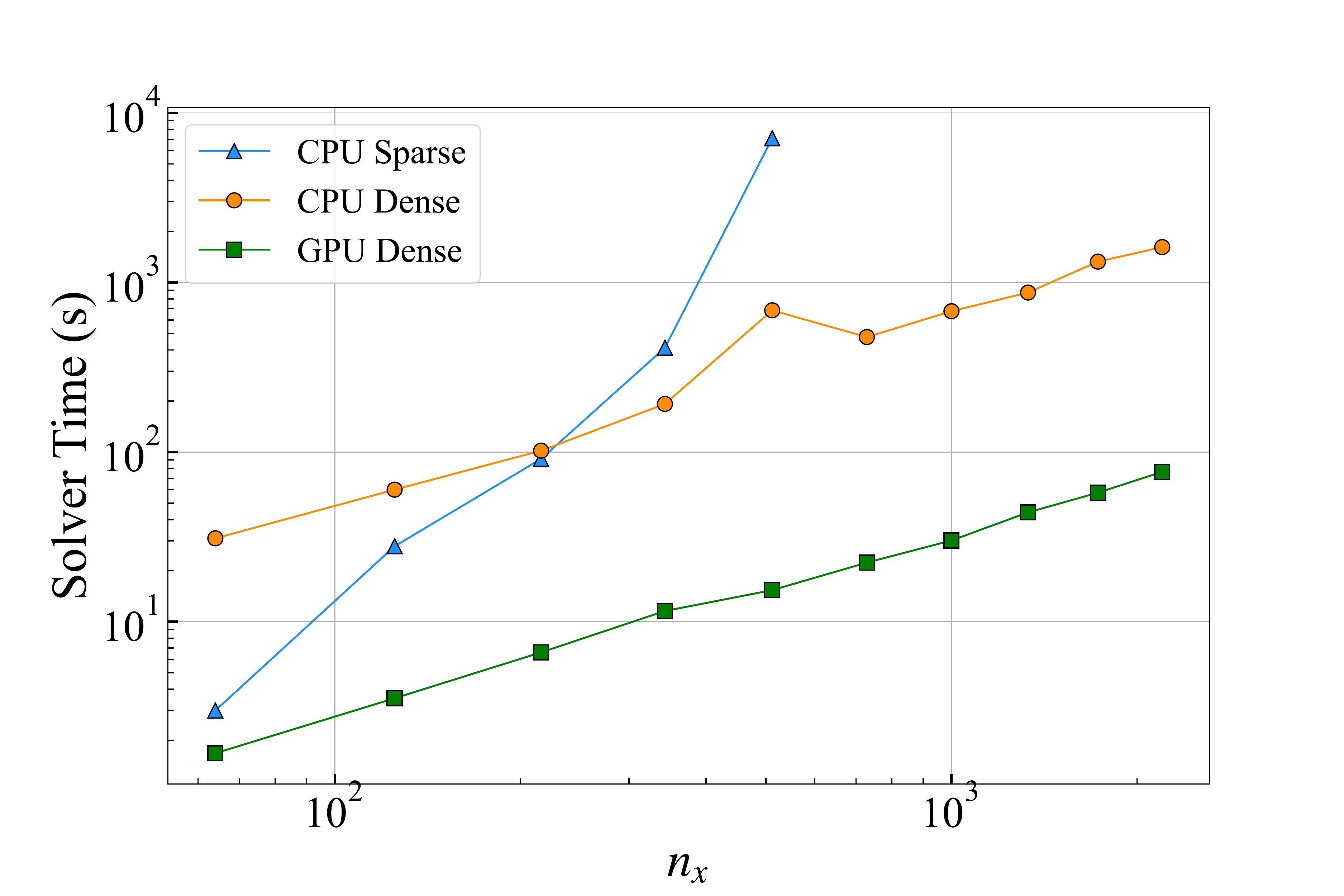}
    \caption{Total solver time for $T = 250$.}
    \label{fig:N250}
\end{figure}

\begin{figure}[t]
    \centering
    \includegraphics[scale=.22]{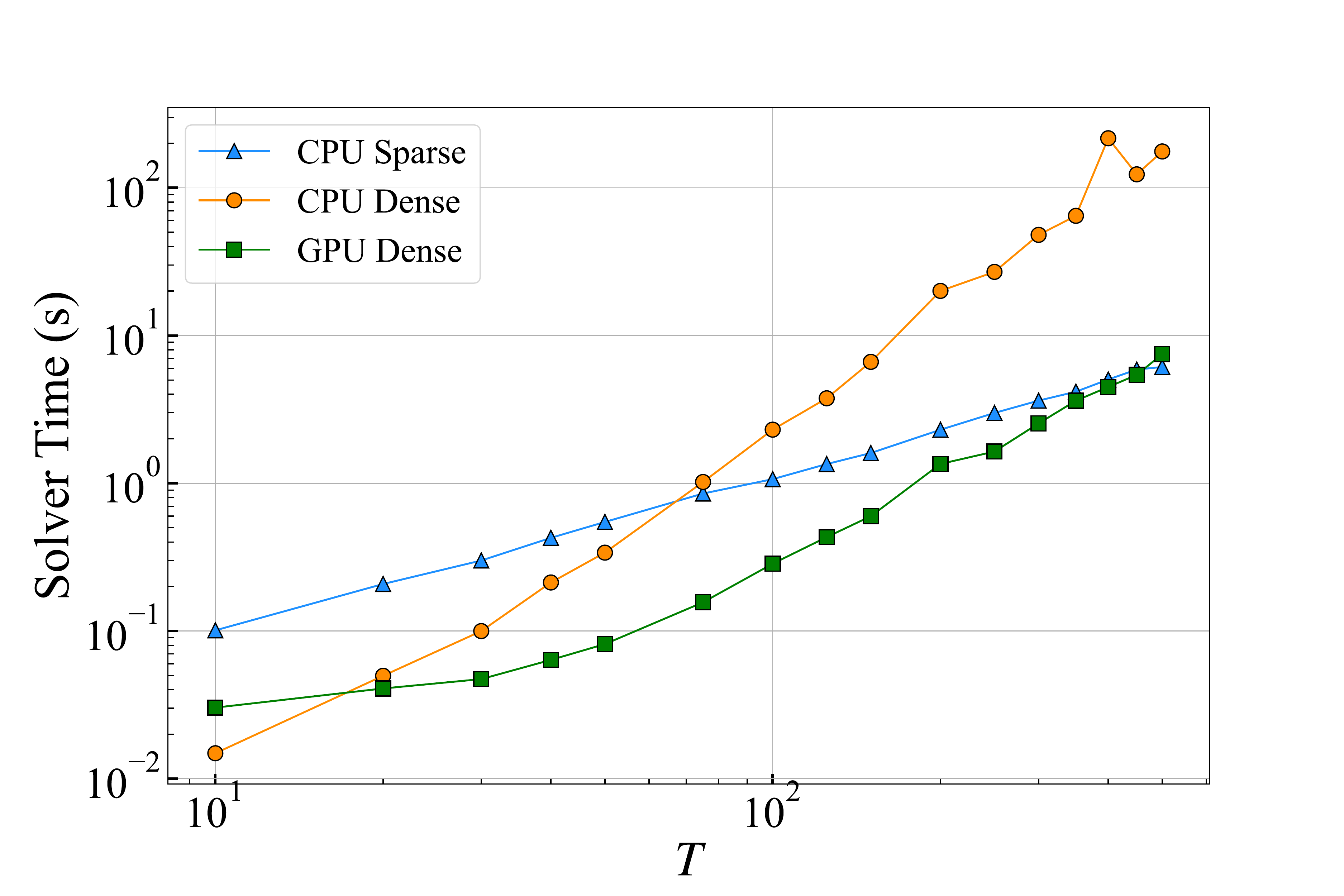}
    \caption{Total solver time for $N = 4$ and varied $T$.}
    \label{fig:nz4}
\end{figure}

We solved problems (i) over a range of $N$ values at three different numbers of time steps: $T = 50$, $T = 150$, and $T = 250$ and (ii) over a range of $T$ values for fixed $N=4$. The results can be seen in Figures \ref{fig:N250} and \ref{fig:nz4} and Tables \ref{tab:N250}--\ref{tab:nz4}. In the tables, ``iter" is the number of iterations within IPMs reported by the solvers ({\tt Ipopt} for the sparse problems and {\tt MadNLP.jl} for the dense problems), ``tot" is the total solver time, and ``lin" is the linear solver time.

\begin{table}[t]
    \caption{Solver statistics for $T=250$ and varied $N$. }
    \setlength{\tabcolsep}{4pt}
    \centering
    \begin{tabular}{c|crr|ccc|ccc}
    \hline\hline
      & \multicolumn{3}{c|}{CPU Sparse} & \multicolumn{3}{c|}{CPU Dense} &  \multicolumn{3}{c}{GPU Dense}\\
    \hline
     $N$ & iter & tot & lin & iter & tot & lin  & iter & tot & lin\\
    \hline
    \hline
    4 & 26 & 3.0  & 2.6  & 28 & 31.0  & 1.10 & 28 & 1.68 & 0.037\\
    5 & 30 & 27.9  & 24.8  & 33 & 60.1  & 1.30 & 33 & 3.54 & 0.045\\
    6 & 29 & 91.2  & 81.4  & 31 & 102.0 & 1.25 & 31 & 6.59 & 0.041\\
    7 & 32 & 412.3 & 372.0 & 38 & 192.3 & 1.55 & 38 & 11.6 & 0.055\\
    8 & 31 & 7086.  & 6617.  & 37 & 685.4 & 1.66 & 37 & 15.4 & 0.051\\
    9 & - & - & -          & 38 & 476.8 & 1.77 & 38 & 22.3 & 0.056\\
    10 & - & - & -         & 37 & 677.2 & 1.76 & 37 & 30.1 & 0.052\\
    11 & - & - & -         & 38 & 872.6 & 1.74 & 38 & 44.1 & 0.058\\
    12 & - & - & -         & 38 & 1328  & 1.74 & 38 & 57.6 & 0.064\\
    13 & - & - & -         & 39 & 1617  & 1.79 & 39 & 76.5 & 0.067\\
    \hline
    \end{tabular}
    \label{tab:N250}
\end{table}

\begin{table}[t]
    \caption{Solver statistics for $T=150$ and varied $N$.}
    \setlength{\tabcolsep}{4pt}
    \centering
    \begin{tabular}{c|crr|ccc|ccc}
    \hline\hline
      & \multicolumn{3}{c|}{CPU Sparse} & \multicolumn{3}{c|}{CPU Dense} &  \multicolumn{3}{c}{GPU Dense}\\
    \hline
     $N$ & iter & tot & lin & iter & tot & lin  & iter & tot & lin\\
    \hline
    \hline
    4 & 25 & 1.7  & 1.4  & 27 & 8.90  & 0.31 & 27 & 0.68 & 0.022\\
    5 & 29 & 17.2  & 15.4  & 32 & 17.4  & 0.54 & 32 & 1.26 & 0.028\\
    6 & 29 & 55.9  & 49.9  & 31 & 28.9  & 0.37 & 31 & 1.93 & 0.027\\
    7 & 30 & 238.2 & 214.1 & 36 & 53.0  & 0.43 & 36 & 4.19 & 0.034\\
    8 & 31 & 4,953.  & 4,636.  & 37 & 155.5 & 0.45 & 37 & 6.12 & 0.034\\
    9 & 31 & 10,357. & 9,703.  & 38 & 110.5 & 0.40 & 38 & 8.05 & 0.034\\
    10 & - & - & -         & 36 & 153.2 & 0.39 & 36 & 10.8 & 0.034\\
    11 & - & - & -         & 38 & 221.1 & 0.43 & 38 & 15.1 & 0.037\\
    12 & - & - & -         & 38 & 353.7 & 0.48 & 38 & 18.9 & 0.040\\
    13 & - & - & -         & 39 & 404.1 & 0.50 & 39 & 25.2 & 0.045\\
    14 & - & - & -         & 40 & 559.6 & 0.53 & 40 & 32.6 & 0.045\\
    \hline
    \end{tabular}
    \label{tab:N150}
  \end{table}

  \begin{table}[t]
    \caption{Solver statistics for $T=50$ and varied $N$. }
    \setlength{\tabcolsep}{4pt}
    \centering
    \begin{tabular}{c|crr|ccc|ccc}
    \hline\hline
      & \multicolumn{3}{c|}{CPU Sparse} & \multicolumn{3}{c|}{CPU Dense} &  \multicolumn{3}{c}{GPU Dense}\\
    \hline
     $N$ & iter & tot & lin & iter & tot & lin  & iter & tot & lin\\
    \hline
    \hline
    4 & 26 & 0.6  & 0.5  & 27 & 0.45 & 0.020 & 27 & 0.09 & 0.008\\
    5 & 28 & 5.0  & 4.4  & 32 & 1.01 & 0.023 & 32 & 0.17 & 0.009\\
    6 & 27 & 17.5  & 15.5  & 29 & 1.66 & 0.022 & 29 & 0.23 & 0.008\\
    7 & 27 & 69.1  & 61.4  & 33 & 3.36 & 0.024 & 33 & 0.41 & 0.010\\
    8 & 32 & 689.0 & 641.1 & 37 & 10.1 & 0.030 & 37 & 0.70 & 0.011\\
    9 & 30 & 3,100.  & 2,936.  & 37 & 7.79 & 0.033 & 37 & 0.94 & 0.011\\
    10 & - & - & -         & 36 & 9.69 & 0.033 & 36 & 1.71 & 0.011\\
    11 & - & - & -         & 37 & 13.9 & 0.034 & 37 & 2.00 & 0.012\\
    12 & - & - & -         & 37 & 18.1 & 0.034 & 37 & 2.43 & 0.012\\
    13 & - & - & -         & 38 & 22.6 & 0.035 & 38 & 3.37 & 0.013\\
    14 & - & - & -         & 38 & 25.6 & 0.033 & 38 & 4.08 & 0.013\\
    \hline
    \end{tabular}
    \label{tab:N50}
\end{table}

\begin{table}[t]
    \caption{Solver statistics for $N=4$ and varied $T$. }
    \setlength{\tabcolsep}{4pt}
    \centering
    \begin{tabular}{c|ccc|ccc|ccc}
    \hline\hline
      & \multicolumn{3}{c|}{CPU Sparse} & \multicolumn{3}{c|}{CPU Dense} &  \multicolumn{3}{c}{GPU Dense}\\
    \hline
     $T$ & iter & tot & lin & iter & tot & lin  & iter & tot & lin\\
    \hline
    \hline
    10  & 26 & 0.10 & 0.079 & 28 & 0.02  & 0.001 & 28 & 0.03 & 0.004\\
    20  & 25 & 0.21 & 0.170 & 27 & 0.05  & 0.002 & 27 & 0.04 & 0.005\\
    30  & 24 & 0.30 & 0.249 & 26 & 0.10  & 0.005 & 26 & 0.05 & 0.005\\
    40  & 25 & 0.43 & 0.356 & 27 & 0.21  & 0.009 & 27 & 0.06 & 0.006\\
    50  & 26 & 0.55 & 0.460 & 27 & 0.34  & 0.016 & 27 & 0.08 & 0.007\\
    75  & 27 & 0.85 & 0.724 & 29 & 1.02  & 0.051 & 29 & 0.16 & 0.010\\
    100 & 25 & 1.07 & 0.900 & 27 & 2.31  & 0.096 & 27 & 0.29 & 0.015\\
    125 & 25 & 1.35 & 1.14  & 27 & 3.76  & 0.162 & 27 & 0.43 & 0.018\\
    150 & 25 & 1.60 & 1.36  & 27 & 6.63  & 0.243 & 27 & 0.60 & 0.021\\
    200 & 26 & 2.30 & 1.97  & 28 & 20.0  & 0.529 & 28 & 1.35 & 0.029\\
    250 & 26 & 2.99 & 2.58  & 28 & 27.0  & 0.979 & 28 & 1.64 & 0.036\\
    300 & 28 & 3.63 & 3.15  & 30 & 48.0  & 1.61  & 30 & 2.54 & 0.047\\
    350 & 27 & 4.17 & 3.61  & 29 & 64.6  & 2.29  & 29 & 3.62 & 0.062\\
    400 & 27 & 5.04 & 4.16  & 29 & 216.4 & 3.30  & 29 & 4.48 & 0.053\\
    450 & 27 & 5.89 & 5.14  & 29 & 123.3 & 4.35  & 29 & 5.41 & 0.067\\
    500 & 27 & 6.12 & 5.30  & 29 & 176.2 & 5.64  & 29 & 7.49 & 0.067\\
    \hline
    \end{tabular}
    \label{tab:nz4}
\end{table}

In all three tests with varying $N$, solution times increased significantly with $N$. However, the results indicate that the computation time of the sparse problem grows more rapidly in $N$ compared with the dense problems. The reason is that in the 3-D temperature control problem, the number of variables in the dense problem remains unchanged (the number of inputs is independent of $N$), whereas the number of variables of the sparse problem scales cubically in $N$. Note that problems with rapidly growing numbers of states can easily arise in discretized PDE-constrained control problems. Also, because of the structure of the discretization of 3-D PDEs, the fill-ins in the factorization grow rapidly, and the sparse solver becomes inefficient as $N$ grows. Although the number of variables in the condensed-space approaches does not change, the solution time is moderately affected by $N$. The reason is that the number of constraints increases with $N$,  making the number of rows in the Jacobian grow linearly with $N$. While the dense system size is still independent of $N$, the condensation procedure, which involves computing $\bJ^\top \bSigma_s\bJ$, is affected by the dimension of $\bJ$.

In all three cases of changing $N$, the dense formulation on the GPU was much faster than the same formulation/algorithm running on the CPU. For $T=250$, the GPU was on average 22 times faster;  for $T=150$, it was 16 times faster; and for $T=50$, it was only 7 times faster. It appears that the GPU operates proportionally faster than the CPU on the larger problem sizes, which is not surprising since the parallelizable computations increase with problem size. However, the solution time for the dense formulation increases faster than the sparse solution time for changing $T$, as can be seen in Figure \ref{fig:nz4} and Table \ref{tab:nz4}. Jerez et al.~\cite{jerez2012sparse} report that the number of flops for computing an interior-point iteration of the dense formulation is $\mathcal{O}\big(T^3 n_u^2 (n_u + n_c)\big)$ while the number of flops for the sparse formulation is $\mathcal{O}\big(T (n_x + n_u)^2 (n_x + n_u + n_c)\big)$. Thus, even though the number of flops for the sparse formulation goes as $\mathcal{O}\big(n_x^3 \big)$,  the sparse formulation can overtake the dense formulation (on either the CPU or GPU) since the dense formulation goes as $\mathcal{O}\big(T^3\big)$ while the sparse formulation goes only as $\mathcal{O}\big(T \big)$. This effect can be seen for $N = 4$ and $T = 500$ (see Figure \ref{fig:nz4} and Table \ref{tab:nz4}).

These results suggest that solving linear-quadratic MPC problems on the GPU using the above methods can reduce the solution times by an order of magnitude for certain cases. This speedup can enable applying MPC to new systems that previously would have been limited by time constraints. The types of problems to which this method can practically be applied are problems with a high ratio of states to inputs and moderate horizon length. As seen above, control systems that can be modeled with PDEs could be one prime candidate. However, our results also suggest that this framework is less effective (and eventually impractical) as the ratio of states to inputs decreases. For example, in Figure \ref{fig:N250} the sparse formulation with $N = 4$ is comparable to the GPU solution time and much faster than the dense formulation on the CPU. The reason is that the problem sizes between the sparse and dense formulations are more similar. Numerous linear-quadratic MPC problems have much lower ratios of states to inputs than that of the temperature control problem presented here, and these problems would be less practical in this framework.

\section{Conclusions and Future Work}

For problems with large numbers of states and few inputs and moderate horizon length, we can speed up linear-quadratic MPC by an order of magnitude by exploiting the structure using GPU/SMID architectures. We introduced the package {\tt DynamicNLPModels.jl}, which reduces the linear-quadratic MPC problem based on user-defined data. The KKT system of this dense model is a reduced KKT system (as compared with the sparse formulation), which {\tt MadNLP.jl} automatically constructs when the model is passed from {\tt DynamicNLPModels.jl} to {\tt MadNLP.jl}. The KKT system can be further condensed within {\tt MadNLP.jl}, resulting in a dense positive definite matrix that can be efficiently factorized on the GPU. Based on the example problem of the heating of a cube, this method can result in significant speedup over the dense form or the sparse form on the CPU.

One of the challenges to be addressed in the future is the GPU memory footprint in building the dense problem. The size of the Jacobian is $(n_x + n_u + n_c) T \times n_u T$. This can become prominently large as $T$ and $n_x$ grow. From the definition of $\bJ$ and $\bB$, however, all the data in $\bJ$ is found in its first $n_u$ columns, making it possible to store the Jacobian implicitly. This could significantly decrease the memory needed for storing $\bJ$. However, building this kernel so that it is competitive with the explicit Jacobian is nontrivial and an area of future work.

% In addition, Wright \cite{wright1991partitioned} showed how the linear MPC problem with time-varying data matrices (e.g., $A$, $B$, $Q$, $R$, etc.) arises from the general discrete-time optimal control problem. The framework shown herein for solving linear MPC problems could potentially be expanded to solve nonlinear optimal control problems. In the future, we would like to see the capabilities of {DynamicNLPModels.jl} expanded to work effectively with general nonlinear problems.


%%%%%%%%%%%%%%%%%%%%%%%%%%%%%%%%%%%%%%%%%%%%%%%%%%%%%%%%%%%%%%%%%%%%%%%%%%%%%%%%

%%%%%%%%%%%%%%%%%%%%%%%%%%%%%%%%%%%%%%%%%%%%%%%%%%%%%%%%%%%%%%%%%%%%%%%%%%%%%%%%

%%%%%%%%%%%%%%%%%%%%%%%%%%%%%%%%%%%%%%%%%%%%%%%%%%%%%%%%%%%%%%%%%%%%%%%%%%%%%%%%

%%%%%%%%%%%%%%%%%%%%%%%%%%%%%%%%%%%%%%%%%%%%%%%%%%%%%%%%%%%%%%%%%%%%%%%%%%%%%%%%

\bibliography{root.bib}
\vspace{0.1cm}
\begin{flushright}
	\scriptsize \framebox{\parbox{2.5in}{Government License: The
			submitted manuscript has been created by UChicago Argonne,
			LLC, Operator of Argonne National Laboratory (``Argonne").
			Argonne, a U.S. Department of Energy Office of Science
			laboratory, is operated under Contract
			No. DE-AC02-06CH11357.  The U.S. Government retains for
			itself, and others acting on its behalf, a paid-up
			nonexclusive, irrevocable worldwide license in said
			article to reproduce, prepare derivative works, distribute
			copies to the public, and perform publicly and display
			publicly, by or on behalf of the Government. The Department of Energy will provide public access to these results of federally sponsored research in accordance with the DOE Public Access Plan. http://energy.gov/downloads/doe-public-access-plan. }}
	\normalsize
\end{flushright}

\end{document}